\documentclass[11pt,twoside]{article}
\usepackage{amssymb}
\usepackage{latexsym}
\usepackage[dvips]{graphicx}
\usepackage{color} 
\newcommand{\blue}[1]{{\color{blue}#1}}

\newcommand{\boldgreek}[1]{\mbox{\boldmath$#1$}}

\newcommand{\R}{I\kern-0.37emR}

\newcommand{\ny}{n\rightarrow\infty}
\newcommand{\Q}{I\kern-0.37emP}
\newcommand{\E}{I\kern-0.37emE}
\newtheorem{THE}{Theorem}[section]

\newtheorem{CO}{Corollary}[section]

\def\bgk{\boldgreek}

\def\bbe{\bgk \beta}
\date{}
\title{Process of the slope components of $\alpha$-regression quantile}
\author{\textsc{Jana Jure\v{c}kov\'a}\thanks{The research was supported by the Grant GA\v{C}R 18-01137S}
\\[2mm]
{\it \normalsize{The Czech Academy of Sciences, Institute of Information Theory and Automation}}\\  \it{\normalsize{Charles University, Faculty of Mathematics and Physics}}\\ \it{\normalsize{Prague, Czech Republic}}}
\begin{document}
\maketitle
\thispagestyle{empty}


\begin{abstract}
We consider the linear regression model along with the process of its $\alpha$-regression quantile, $0<\alpha<1$. We are interested mainly in the slope components of $\alpha$-regression quantile and in their dependence on the choice of $\alpha.$ 
While they are invariant to the location, and only  
the intercept part of the $\alpha$-regression quantile estimates the quantile $F^{-1}(\alpha)$ of the model errors, their dispersion depends on $\alpha$ and is infinitely increasing as $\alpha\rightarrow 0,1$, in the same rate as for the ordinary quantiles.
 
We study the  process of $R$-estimators of the slope parameters over $\alpha\in[0,1]$, generated by the H\' ajek  rank scores.   
We show that this process, standardized by $f(F ^{-1}(\alpha))$ under exponentially tailed $F$,
  converges to the vector of independent Brownian bridges. The same course is true for the process of the slope components of $\alpha$-regression quantile. \\


\end{abstract}


\section{Introduction}
\setcounter{equation}{0}

We start with the linear regression model
\begin{equation}
\label{1}
Y_{ni}=\beta_0+{\bf x}_{ni}^{\top}{\boldgreek\beta}+e_{ni}, \quad
i=1,\ldots,n
\end{equation}
with observations $Y_{n1},\ldots,Y_{nn},$
independent errors $e_{n1},\ldots,e_{nn},$
identically distributed according to an unknown distribution function $F;$
with the vector of covariates\\ ${\bf x}_{ni}=(x_{i1},\ldots,x_{in})^{\top},$ 
$i=1,\ldots,n,$  unknown parameter
${\boldgreek\beta}=(\beta_1,\ldots,\beta_p)^{\top}$ of interest, and nuisance intercept $\beta_0.$ 
The regression $\alpha$-quantiles, $0\leq\alpha\leq 1,$ introduced in \cite{KB78}, are an important tool mainly in economics, where the quantile regression became a technical term. 
Remind that the regression $\alpha$-quantile of model (\ref{1}) is defined as the solution of the minimization
\begin{eqnarray}\label{RQ}
&&(\hat{\beta}_0(\alpha), \hat{\boldgreek\beta}(\alpha)=\\
&&=\arg\min\{\alpha\sum_{i=1}^n(Y_i-b_0-\mathbf x_i^{\top}\mathbf b)^+
+(1-\alpha)\sum_{i=1}^n(Y_i-b_0-\mathbf x_i^{\top}\mathbf b)^-\},\nonumber\\
&&\qquad\qquad\qquad\qquad b_0\in\mathbb R_1, \; \mathbf b\in\mathbb R_p, \;  0<\alpha< 1.\nonumber
\end{eqnarray}
The population counterpart of (\ref{RQ}) is 
$$\{\tilde{\beta_0}(\alpha)=\beta_0+F^{-1}(\alpha), \; \widetilde{\boldgreek\beta}(\alpha)=(\beta_1,\ldots,\beta_p)^{\top}, \; 0<\alpha<1\}.$$
Hence,  only the intercept part of the $\alpha$-regression quantile reflects the quantile $F^{-1}(\alpha)$ of the probability distribution $F,$ while $\widetilde{\boldgreek\beta}(\alpha)$ only reflects the slopes. If $\alpha$ runs over the interval $(0,1),$ we get the regression quantile process with step-functions trajectories with  number of breakpoints increasing with the number $n$ of observations. There is a rich literature devoted to the concepts connected with regression quantile, its processes and applications. As an excellent review we recommend Koenker's book \cite{Koenkerbook}.  
The choice of $\alpha$ is an important decision, namely when $Y_i$ reflects the  loss, when we should consider the balance between the underestimation and overestimation of our risk.  For the applications is important the shape of the limiting processes over $0<\alpha<1$ and the shape of various functionals of the regression quantile, characterizing the economic properties. 

Alternatively to the regression quantile, we can follow  the intercept and slope components separately. The so called two-step regression quantile, proposed in \cite{twostep}, first estimates the slope components $\boldgreek\beta$ with the aid of rank estimator $\widehat{\boldgreek\beta}_{nR}$ and then estimates the intercept as the $\alpha$-quantile of the residuals $Y_i-\mathbf x_i^{\top}\widehat{\boldgreek\beta}_{nR},$ $i=1,\ldots,n.$ The two-step regression quantile process is asymptotically equivalent to the ordinary regression quantile process, only its number of breakpoints  differs, being exactly $n.$ The empirical processes corresponding to the regression quantiles and their inversions are numerically illustrated in \cite{aplimat}.

The R-estimate of slopes $\boldgreek\beta$
 is generally defined as the minimizer  $\widehat{\bbe}_{nR}$ of
the Jaeckel \cite{Jaeckel1972} measure of the rank dispersion
\begin{eqnarray}\label{11d}
\mathcal D_{n}(\mathbf b)&=&\sum_{i=1}^n(Y_{ni}-\mathbf x_{ni}^{\top}\mathbf b) \left(a_n(R_{ni}(Y_i-\mathbf x_{ni}^{\top}\mathbf b)-\bar{a}_n\right)\\
&=&\sum_{i=1}^n\left[(Y_{ni}-\bar{Y}_n)-(\mathbf x_{ni}-\bar{\mathbf x}_n)^{\top}\mathbf b)\right]a_n(R_{ni}(Y_i-\mathbf x_{ni}^{\top}\mathbf b), \; \mathbf b\in\mathbb{R}_p,\nonumber
\end{eqnarray}
where
$ R_{ni}(Y_{ni}-\mathbf x_{ni}^{\top}\mathbf b)$  is the rank of the residual 
$Y_{ni}-\mathbf x_{ni}^{\top}\mathbf b $ among $Y_{n1}-\mathbf x_{n1}^{\top}\mathbf b,%
\ldots,Y_{nn}-\mathbf x_{nn}^{\top}\mathbf b,$  
$\bar{Y}_n=\frac 1n\sum_{i=1}^n Y_{ni}$, $\bar{\mathbf x}_n
=\frac 1n\sum_{i=1}^n\mathbf x_{ni},$ $a_n(i)$ are the scores and $\bar{a}_n=\frac 1n\sum_{i=1}^na_n(i)$.
Notice that $\widehat{\boldgreek\beta}_{nR}$ in invariant to the shift in location, hence it is independent of $\beta_0.$

The scores $a_n(i), \; i=1,\ldots,n$ are typically generated by a
function $\varphi(u):(0,1)\mapsto \mathbb R_1,$ nondecreasing and square integrable on (0,1), such that
  $$\lim_{\ny}\int_0^1\Big(a_n(1+[nu]))-\varphi(u)\Big)^2 du =0.$$  
For instance, $a_n(i)=\varphi\left(\frac{i}{n+1}\right), \; i=1,\ldots,n.$ 

Particularly, we shall consider the following family of score functions $\Big\{\varphi_{\alpha}(u), \; 0\leq\alpha\leq 1, \; 0\leq u \leq 1\Big\}:$
$$\varphi_\alpha(u)=\left\{\begin{array}{lll}
                 0 & \ldots & 0\leq u\leq\alpha\leq 1\\
                 1 & \ldots & 1\geq u>\alpha \geq 0.\\
                 \end{array} \right. 
$$ 
 As $\ny,$ the function $\varphi_\alpha(u)$ generates the following scores :    
\begin{equation}
\label{8b}
a_n(i,\alpha)=\left\{\begin{array}{lll}
              0 & \ldots & i\leq n\alpha\\
              i-n\alpha & \ldots & n\alpha \leq i \leq n\alpha +1\\
              1 & \ldots & i \geq n\alpha +1\\
              \end{array}
              \right.                            
\end{equation}
$i=1,\ldots,n.$ Notice that $a_n(i,\alpha)$ is continuous in $\alpha\in(0,1).$ The scores $a_n(i,\alpha), \; i=1,\ldots,n$ are known as H\' ajek's rank scores  (see H\' ajek \cite{Hajek} and 
H\' ajek and \v{S}id\' ak \cite{HS1967}). 

If $R_{n1},\ldots,R_{nn}$ are the ranks of  random variables $Z_{1},\ldots,Z_{n},$ then the vector\\ $(a(R_{n1},\alpha),\ldots,a(R_{nn},\alpha))$ is a solution of the linear programming
\begin{eqnarray}\label{programming}
&&\sum_{i=1}^n Z_i~{a}_n(R_{ni},\alpha)=\max\\
\mbox{under } &&\sum_{i=1}^n {a}_n(R_{ni},\alpha)=n(1-\alpha)\\
&&0\leq {a}_n(R_{ni},\alpha)\leq 1, \; \; i=1,\ldots,n\nonumber
\end{eqnarray}
(cf. also \cite{GJ,GJKP}).
In this case the Jaeckel criterion (\ref{11d}) asymptotically simplifies, as $\ny$, to 
\begin{eqnarray}\label{11ddd}
&&\mathcal D_{n\alpha}(\mathbf b)=\nonumber\\ 
&&=\sum_{i=1}^n\left[(Y_{ni}-\bar{Y}_n)-(\mathbf x_{ni}-\bar{\mathbf x}_n)^{\top}\mathbf b)\right]\Big(I[R_{ni}(Y_i-\mathbf x_{ni}^{\top}\mathbf b)\geq n\alpha]
\nonumber\\
&&+(R_{ni}(Y_i-\mathbf x_{ni}^{\top}\mathbf b)-n\alpha)I[n\alpha \leq R_{ni}(Y_i-\mathbf x_{ni}^{\top}\mathbf b) \leq n\alpha +1]\Big)\\
&&\approx\sum_{i=1}^n\left[(Y_{ni}-\bar{Y}_n)-(\mathbf x_{ni}-\bar{\mathbf x}_n)^{\top}\mathbf b)\right]I[R_{ni}(Y_i-\mathbf x_{ni}^{\top}\mathbf b)\geq n\alpha].\nonumber
\end{eqnarray}
{Jaeckel proved that ${\cal D}_{n\alpha}({\mathbf b})$ is continuous, convex and piecewise linear function of ${\mathbf b}\in \mathbb R_p,$ thus differentiable with gradient
\begin{eqnarray}\label{301}
	&&\frac{\partial{\mathcal D}_{n\alpha}({\mathbf b})}{\partial{\mathbf b}}\Big|_{\mathbf b_0}
	=	-\sum_{i=1}^n(\mathbf x_{ni}-\bar{\mathbf x}_n)I\left[R_{ni}(Y_i-\mathbf x_{ni}^{\top}\mathbf b_0)\geq n\alpha \right] 
	  \end{eqnarray}
at any point ${\mathbf b}_0\in\mathbb R_p$ of differentiability.} 
Notice that the gradients of the Jaeckel measure are just the H\' ajek scores.  Using the uniform asymptotic linearity of the H\' ajek scores, (see Theorem \ref{Theorem11}), we can approximate the Jaeckel measure by a quadratic function. 
 
Our subject of interest is to investigate the possible convergence of the process of R-estimators $\{\widehat{\boldgreek\beta}_{n\alpha}, \; 0<\alpha<1\}$, generated by the H\'ajek scores (\ref{8b}). H\' ajek and \v{S}id\' ak \cite{HS1967}  proved the weak convergence of the process of H\' ajek's scores to the Brownian bridge, under the i.i.d. observations as well as under contiguous (Pitman) alternatives. The intercept component of the $\alpha$-regression quantile reflects the population quantile $F^{-1}(\alpha),$ but it is not the case of the slope components. However, the dispersion of the process of slopes expands for $\alpha\rightarrow 0,1$ and its variance copies the variance of the $\alpha$-quantile, i.e. $\alpha(1-\alpha)f^{-2}(F^{-1}(\alpha)).$ Under conditions  on the tails of distribution of model errors, such as
imposed in \cite{GJKP}, we can prove the weak convergence of the  process  $\{f(F^{-1}(\alpha))(\widehat{\bbe}_{n\alpha}-\boldgreek\beta)\}$ to the vector of independent Brownian bridges over 
the compact subsets of $[0,1].$

\section{Process of R-estimates of slopes and its asymptotics}
\setcounter{equation}{0}

Let $\widehat{\boldgreek\beta}_{n\alpha}$ be the R-estimator of $\boldgreek\beta,$ based on the H\' ajek rank scores, i.e. the minimizer of (\ref{11ddd}). Following the steps of \cite{GJKP}, we shall first  study the order of $\widehat{\boldgreek\beta}_{n\alpha}$ over $(\alpha_n^*, 1-\alpha_n^*)$ and show that the process of R-estimators converges to the vector of independent Brownian bridges for some $\alpha_n^*\downarrow 0$ as $\ny.$ This, in turn, will lead to the convergence over $\alpha\in(\alpha_0,1-\alpha_0)$ with any $0<\alpha_0<1/2$ fixed.

Consider the  process of the H\' ajek rank scores
\begin{eqnarray}\label{J2.1a}
&&\mathcal A_n(n^{-1/2}\mathbf b)\\
&&=\left\{\mathcal A_{n\alpha}(n^{-1/2}\mathbf b)=n^{-1/2}\sum_{i=1}^n (\mathbf x_{ni}-\bar{\mathbf x}_n)a_{n\alpha}(R_{ni}(Y_i-n^{-1/2}\mathbf x_{ni}^{\top}\mathbf b),\alpha): \;  \; 0\leq \alpha\leq 1\right\}\nonumber
\end{eqnarray}
for $\mathbf b\in\mathbb R_p.$ The R-estimator $\widehat{\bbe}_{n\alpha}$ is the minimizer of $\mathcal D_{n\alpha}(n^{-1/2}\mathbf b)$ and $\mathcal A_{n\alpha}(n^{-1/2}\mathbf b)$ is its gradient, due to (\ref{11ddd}) and (\ref{301}).  The results in \cite{GJKP} and \cite{Jur1992a} imply that the process (\ref{J2.1a}) is uniformly asymptotically linear in $\mathbf b,$ what enables to approximate $\mathcal D_{n\alpha}(n^{-1/2}\mathbf b)$ by a quadratic function and then to approximate $\widehat{\bbe}_{n\alpha}$ by its minimizer.

In order to realize these approximations, we impose the following  conditions on the distribution of the model errors and on the triangular array  of covariates ${\mathbf x}_{n1},\ldots,{\mathbf x}_{nn}$. These conditions are only sufficient and apparently can be weakened.
\begin{description}
\item{(F1)} The density $f(x)=F^{\prime}(x)$ is absolutely continuous and bounded with bounded derivative $f^{\prime}$ for $A<x<B,$ where $-\infty\leq A=%
	\sup\{x: F(x)=0\}$ and 
	$+\infty\geq B=\inf\{x: F(x)=1\}.$ 
	\item{(F2)} The density $f(x)=F^{\prime}(x)$ 
	 is monotonically decreasing as $x\downarrow A$ or $x\uparrow B$ and $f^{\prime}(x)/f(x)|\leq c|x|$ for $x\geq K(\geq 0), \; c>0.$
	\item{{(F3)}}  $|F^{-1}(\alpha)|\leq c(\alpha(1-\alpha))^{-a}$ and 
	similarly, $1/f(F^{-1}(\alpha))\leq c(\alpha(1-\alpha))^{-a-1}$ for $0<\alpha\leq \alpha_0$ and $1-\alpha_0\leq\alpha<1$ where 
	$0<a<\frac  14-\varepsilon,$ $\varepsilon>0, \;  0<\alpha_0\leq 1/2.$	
\end{description}
\begin{description}
	\item{(X1)} 
The matrix 
$$\mathbf Q_n=\sum_{i=1}^n({\mathbf x}_{ni}-\bar{\mathbf x}_n)({\mathbf x}_{ni}-\bar{\mathbf x}_n)^{\top}, \quad \bar{\mathbf x}_n=\frac 1n\sum_{i=1}^n{\mathbf x}_{ni}$$
has the rank $p$ and $n^{-1}\mathbf Q_n\rightarrow \mathbf C$ as $\ny,$
where $\mathbf C$ is a positively definite $p\times p$ matrix. Moreover, we assume
\begin{equation}\label{Noether}
\lim_{\ny}\max_{1\leq i\leq n}({\mathbf x}_{ni}-\bar{\mathbf x}_n)^{\top}\mathbf Q_n^{-1}({\mathbf x}_{ni}-\bar{\mathbf x}_n)=0 \quad\mbox{ (\textit{Noether condition}).}
\end{equation}
	\item{(X2)} $n^{-1}\sum_{i=1}^n\|\mathbf x_{ni}\|^4 =O(1)$ as $\ny,$ and\\
	$\max_{1\leq i\leq n}\|\mathbf x_{ni}\|=O\left(n^{(2(b-a)-\delta)/(1+4b)}\right)$ as $\ny$ for some $b>0, \; \delta>0$ such that
	$0<b-a<\frac{\varepsilon}{2}.$
\end{description}
As a consequence of Section V.3.5 in \cite{HS1967},  we get the following weak convergence in the Prokhorov topology under  $\mathbf b=\mathbf 0$ 
\begin{equation}\label{3.5.2a}
\Big\{n^{1/2}\mathbf Q_n^{-1/2}\mathcal A_{n\alpha}(\mathbf 0): 0\leq\alpha\leq 1\Big\}\stackrel{\mathcal D}{\rightarrow}\mathbf W^*_p
\end{equation}
as $\ny,$ where $\mathbf W^*_p$ is the vector of $p$ independent Brownian bridges (see \cite{HS1967} and \cite{GJ}). Furthermore, under a sequence of contiguous alternatives, 
when $Y_{ni}=Y_{ni}^0+n^{-1/2}\mathbf x_{ni}^{\top}\mathbf b,$\\ $i=1,\ldots,n$ with $Y_{ni}^0$ independent having distribution function $F,$ 
there  also applies the following convergence to the vector of $p$ independent Brownian bridges
\begin{eqnarray}\label{VI.3.2.2a}
&&\Big\{n^{1/2}\mathbf Q_n^{-1/2}\mathcal A_n(\alpha,\blue{n^{-1/2}}\mathbf b)-n^{-1/2}\mathbf Q_n^{1/2}\mathbf b f(F^{-1}(\alpha)) \; :  \; 0\leq\alpha\leq 1\Big\}\stackrel{\mathcal D}{\rightarrow}\mathbf W^*_p 
\end{eqnarray}
{ as } $\ny$ (see \cite{HS1967}, Theorem VI.3.2).
   The first result is the uniform asymptotic linearity of $A_n(\alpha,\blue{n^{-1/2}}\mathbf b)$
in $\mathbf b$, proven in \cite{Jur1992a}.

Denote 
\begin{equation}\label{J2.5}
\sigma_{\alpha}=\frac{(\alpha(1-\alpha))^{1/2}}{f(F^{-1}(\alpha))}, \; 0<\alpha<1 \quad\mbox{ and }\quad 
\alpha_n^*=1/n^{1+4b} \quad \mbox{ with} \;\; b \;\; \mbox{from} \;\; (X2).
\end{equation}
\begin{THE}\label{Theorem11} Assume that $F$ and $\mathbf X_n$ satisfy (F1)--(F3) and (X1)--(X2). Then
\begin{eqnarray}\label{J2.2}
&&\sup\Big\{(\alpha(1-\alpha))^{-1/2}\Big |\mathcal A_n(\alpha,\blue{n^{-1/2}\sigma_{\alpha}}\mathbf b)-\mathcal A_n(\alpha,\mathbf 0)+n^{-1}\mathbf Q_n \mathbf b\Big |: \\
&&\qquad\qquad\qquad\qquad \;\; \|\mathbf b\|\leq K, \; \alpha_n^*\leq\alpha\leq 1-\alpha_n^*\Big \}\stackrel{p}{\rightarrow}0\nonumber
\end{eqnarray}
and
\begin{eqnarray}\label{J2.28}
&&\sup\Big\{\Big |\mathcal A_n(\alpha,\blue{n^{-1/2}}\mathbf b)-\mathcal A_n(\alpha,\mathbf 0) 
+f(F^{-1}(\alpha))n^{-1}\mathbf Q_n \mathbf b\Big |:\\
&&\qquad\qquad\qquad\qquad \;\;  \|\mathbf b\|\leq K, \;  0\leq\alpha\leq 1\Big \}\stackrel{p}{\rightarrow}0\nonumber
\end{eqnarray}
as $\ny,$ for any fixed $K, \; 0<K<\infty.$
\end{THE}
\textbf{Proof}. The theorem is proven in \cite{Jur1992a}.\\

The following theorem gives the asymptotic behavior of the R-estimator of slope parameter over the interval $[\alpha_n^*, 1-\alpha_n^*]$.
\begin{THE}\label{LemmaJ1}
Under the conditions of Theorem \ref{Theorem11}, as  $\ny,$
\begin{equation}\label{JJ2.12}
\sup\left\{n^{1/2}\sigma_{\alpha}^{-1}\|\widehat{\boldgreek\beta}_{n\alpha}-\boldgreek\beta\|: \; \alpha_n^*\leq\alpha\leq 1-\alpha_n^*\right\}=O_p(1)
\end{equation}
and
\begin{eqnarray}\label{J2.12}
&&\sup\Big\{\|n^{1/2}\sigma_{\alpha}^{-1}[\widehat{\boldgreek\beta}_{n\alpha}-\boldgreek\beta]-(\alpha(1-\alpha))^{-1/2} \; n\mathbf Q_n^{-1}\mathcal A_{n\alpha}
(\mathbf 0)\|:\\
&&\qquad\qquad\qquad\qquad \;\; \alpha_n^*\leq\alpha\leq 1-\alpha_n^*\Big\}
=o_p(1).\nonumber 
\end{eqnarray}
Moreover, the process
\begin{equation}\label{process}
\left\{f(F^{-1}(\alpha))\mathbf Q_n^{1/2}(\widehat{\boldgreek\beta}_{n\alpha}-\boldgreek\beta): \; \alpha_n^*\leq\alpha\leq 1-\alpha_n^*\right\}
\end{equation}
converges to the vector of independent Brownian bridges. 
\end{THE}
\textbf{Proof.} The theorem is proven in Section 3.\\

As a consequence, we conclude that the process $\{f(F^{-1}(\alpha))\mathbf Q_n^{1/2}(\widehat{\boldgreek\beta}_{n\alpha}-\boldgreek\beta)\}$ converges to the vector of Brownian bridges over the interval $[\alpha_0, 1-\alpha_0]$ for any fixed $0<\alpha_0<1/2$, i.e. over the compact subsets of (0,1).
\begin{CO}\label{corollary}
Under the conditions of Theorem \ref{LemmaJ1}, the process
\begin{equation}\label{process1}
\left\{f(F^{-1}(\alpha))\mathbf Q_n^{1/2}(\widehat{\boldgreek\beta}_{n\alpha}-\boldgreek\beta): \; 0<\alpha< 1\right\}
\end{equation}
converges to the vector of independent Brownian bridges in $\mathcal D(0,1)^p$. The convergence over $0<\alpha<1$ is in the sense that the process converges over the interval $[\alpha_0, 1-\alpha_0]$ for any fixed $0<\alpha_0<1/2$, i.e. converges over the compact subsets of (0,1).
\end{CO}
\section{Proofs}
\setcounter{equation}{0}
\subsection*{Proof of Theorem \ref{LemmaJ1}}
Notice that
\begin{equation}\label{J2.15}
\mathbf t_{n\alpha}= n^{1/2}\sigma_{\alpha}^{-1}(\widetilde{\boldgreek\beta}_{n\alpha}-\boldgreek\beta)
\end{equation}
 minimizes $[\mathcal D_{n\alpha}(n^{-1/2}\sigma_{\alpha}\mathbf b)-\mathcal D_{n\alpha}(\mathbf 0)].$ 
Theorem \ref{Theorem11} leads to the following quadratic approximation of $D_{n\alpha}(\mathbf b):$
\begin{eqnarray}\label{3.1}
&&\sup\Big\{\Big|(\alpha(1-\alpha))^{-1/2}\Big(\sigma_{\alpha}^{-1}[\mathcal D_{n\alpha}(n^{-1/2}\sigma_{\alpha}\mathbf b)-\mathcal D_n(\mathbf 0)]+\mathbf b^{\top}\mathcal A_{n\alpha}(\mathbf 0)\Big)
-\frac 12 n^{-1}\mathbf b^{\top}\mathbf Q_n\mathbf b\Big|:\nonumber\\ 
&&\quad \|\mathbf b\|\leq K, \alpha_n^*\leq \alpha\leq 1-\alpha_n^*\Big\} \stackrel{p}{\rightarrow}0 \; \mbox{ as }\ny \; \mbox{ for any fixed } K>0.
\end{eqnarray}
This further implies that as $\ny$
\begin{eqnarray}\label{J2.16}
&&\min_{\|\mathbf b\|\leq K}\left[(\alpha(1-\alpha))^{-1/2}\sigma_{\alpha}^{-1}[\mathcal D_{n\alpha}(n^{-1/2}\sigma_{\alpha}\mathbf b)-\mathcal D_{n\alpha}(\mathbf 0)]\right]\\
&&=\min_{\|\mathbf b\|\leq K}\Big[\frac 12 n^{-1}\mathbf b^{\top}\mathbf Q_n\mathbf b-(\alpha(1-\alpha))^{-1/2}\mathbf b^{\top}\mathcal A_{n\alpha}(\mathbf 0)\Big]+o_p(1)\nonumber
\end{eqnarray}
 uniformly for $\alpha_n^*\leq\alpha\leq 1-\alpha_n^*,$ for any $K, \; 0<K<\infty.$  
 
Moreover,
\begin{eqnarray}\label{J2.19a}
&&\min_{\mathbf b\in\mathbb R_p}\Big[\frac 12 n^{-1}\mathbf b^{\top}\mathbf Q_n\mathbf b-(\alpha(1-\alpha))^{-1/2}\mathbf b^{\top}\mathcal A_{n\alpha}(\mathbf 0)\Big]\nonumber\\
&&=-\frac 12(\alpha(1-\alpha))^{-1}\mathcal{A}_{n\alpha}^{\top}(\mathbf 0) \; n\mathbf Q_n^{-1}\mathcal{A}_{n\alpha}(\mathbf 0)
\end{eqnarray}
and
\begin{eqnarray}\label{J2.19}
&&\arg\min_{\mathbf b\in\mathbb R_p}\Big[\frac 12 n^{-1}\mathbf b^{\top}\mathbf Q_n\mathbf b-(\alpha(1-\alpha))^{-1/2}\mathbf b^{\top}\mathcal A_{n\alpha}(\mathbf 0)\Big]\nonumber\\
&&=(\alpha(1-\alpha))^{-1/2} \; n\mathbf Q_n^{-1}\mathcal A_{n\alpha}(\mathbf 0)\\
&&=\mathbf u_{n\alpha} \quad \mbox{ (SAY) }.\nonumber 
 \end{eqnarray}
Notice that $\|\mathbf u_{n\alpha}\|=O_p(1)$ uniformly in  $\alpha_n^*\leq\alpha\leq 1-\alpha_n^*$ by (\ref{J2.1a}). Inserting $\mathbf b=\mathbf u_{n\alpha}$ in
 (\ref{3.1}), we obtain
\begin{eqnarray}\label{J2.18}
&&\sup\Big\{\Big|(\alpha(1-\alpha))^{-1/2}\sigma_{\alpha}^{-1}[\mathcal D_{n\alpha}(n^{-1/2}\sigma_{\alpha}\mathbf u_{n\alpha})-\mathcal D_n(\mathbf 0)]\\
 &&+\frac 12 (\alpha(1-\alpha))^{-1}\mathcal{A}_{n\alpha}^{\top}(\mathbf 0) \; n\mathbf Q_n^{-1}\mathcal{A}_{n\alpha}(\mathbf 0)\Big|:
\; \alpha_n^*\leq\alpha\leq 1-\alpha_n^*\Big\}=o_p(1)\nonumber
\end{eqnarray}
Hence, using the convexity of $\mathcal D_n,$ we apply the approach of Pollard in \cite{Pollard} and conclude
\begin{equation}\label{J2.20}
\sup\left\{\|\mathbf t_{n\alpha}-\mathbf u_{n\alpha}\|: \; \alpha_n^*\leq\alpha\leq 1-\alpha_n^*\right\}=o_p(1).
\end{equation}
The convergence of (\ref{process}) to the vector of Brownian bridges follows from (\ref{3.5.2a}).
 \hfill $\Box$\\
\vspace{5mm}

If $1-\alpha\geq 1-\alpha_n^*$, then $R_{ni}(Y_i-\mathbf x_i^{\top}\mathbf b)\geq n(1-\alpha)$ only for the maximal residual $Y_i-\mathbf x_i^{\top}\mathbf b.$ Hence the estimator $\widetilde{\boldgreek\beta}_{n(1-\alpha)}$ minimizes the maximal residual over $\mathbf b\in\mathbb R_p.$ More precisely,
\begin{equation}\label{31}
\widetilde{\boldgreek\beta}_{n(1-\alpha)}=\arg\min_{\mathbf b\in\mathbb R_p}\left\{[Y_{ni}-\bar{Y}_n-(\mathbf x_{ni}-\bar{\mathbf x}_n)^{\top}\mathbf b]_{n:n}\right\}.
\end{equation}
Denote as $D_n$ the antirank of the maximal residual. Then 
$$\left[Y_{nD_n}-\bar{Y}_n-(\mathbf x_{nD_n}-\bar{\mathbf x}_n)^{\top}\widetilde{\boldgreek\beta}_{n(1-\alpha)}\right]\leq [Y_{nD_n}-\bar{Y}_n-(\mathbf x_{nD_n}-\bar{\mathbf x}_n)^{\top}\mathbf b],$$
hence 
\begin{equation}\label{32}
(\mathbf x_{nD_n}-\bar{\mathbf x}_n)^{\top}\widetilde{\boldgreek\beta}_{n(1-\alpha)}\geq (\mathbf x_{nD_n}-\bar{\mathbf x}_n)^{\top}\mathbf b 
\end{equation}
for any $\mathbf b\in\mathbb R_p,$ including $\mathbf 0$ and any other estimator of $\boldgreek\beta.$ Moreover, notice that $\widetilde{\boldgreek\beta}_{n(1-\alpha)}$ is constant for $0<\alpha\leq\alpha_n^*,$ i.e.
$$\widetilde{\boldgreek\beta}_{n(1-\alpha)}=\widetilde{\boldgreek\beta}_{n(1-\alpha_n^*)} \; \mbox{ for } \;  0<\alpha\leq\alpha_n^*.$$ 
Analogously, $\widetilde{\boldgreek\beta}_{n\alpha_0}=\widetilde{\boldgreek\beta}_{n\alpha_n^*}$ for $0<\alpha_0\leq\alpha_n^*$,
hence the convergence holds over $[\alpha_0, 1-\alpha_0]$ for any $0<\alpha_0,1/2$ and thus for the compact subintervals of $(0,1)$.
\hfill $\Box$\\

\end{document}